# Extending iLQR method with control delay

Cheng Ju, Yan Qin and Chunjiang Fu[*]

1. Introduction

In the framework of optimal control, the research of human motor control system has been progressing for 30 years. Numerous experimental phenomena have risen up fresh questions and contributed to the development of theoretical study. It has achieved the range of problems including nonlinearity, redundancy, signal dependant noise and constraint. In last century, feedforward optimal control was introduced first with various criterion of cost (Uno, 1989). Later, the importance of the involvement of noise was indicated by (Harris and Wolpert, 1998). After 20[th] century, feedback was added in to form a relatively unified framework to computationally describe the motor control system, which is called optimal feedback control framework also including multiplicative noise and constraints (Todorov, 2002).

However, such framework does not express how to deal with the time delay. Physiologically speaking, delays exist in all stages of motor system such as receiving sensory data, transferring motion command and muscle responding, which vary from 10ms to 150ms depending on the task. Delays can produce error and even instability of present state estimation, which cannot be neglected in the stochastic optimal feedback control framework (OFC) if the time of motion is long enough. Some literatures even employ delay to explain famous Fitts Law (Beamish, 2006). This work imports time delay into the OFC. As it is known, it is difficult to find a globally-optimal control regulator in such complex problem. By the idea of iterative linear quadratic regulator(iLQR, Weiwei Li, 2005), a discrete-time linear quadratic stochastic control system is developed to approximate the original problem. Then the local optimal feedback control regulator is obtained in a recursive form. This work demonstrates OFC is capable of including time delay in theory.

The paper is organized as follows. Section 2 presents the nonlinear stochastic dynamic system with time delay. An approximate linear quadratic stochastic control system is introduced by linearization techniques. In Section 3, the iterative linear quadratic regulator (iLQR) with input delay algorithm is designed in our main theorem.

2. Problem Formulation

2.1 The non-linear stochastic dynamic system with time-delay

Consider a class of non-linear dynamic system with time-delay described by the stochastic state equation

$$dx(t) = f(x(t), u(t), u(t-\tau))dt + F(x(t), u(t), u(t-\tau))d\omega(t) \tag{1}$$

with the initial condition $x(0) = x_0$ and $u(s) = 0$ for any $s \in [-\tau, 0)$. Here $\tau$ is a positive fixed delay, $x \in R^n$ is the state variable and $u \in R^d$ is the control input variable without any constraint in our paper. $\omega(t)$ is a $d$-dimension standard Brownian motion. The coefficients $f$ and $F$ are continuous with respect to all their arguments. The derivatives of these are continuous and bounded. The cost function to be minimized is defined as follow.

$$J(t,x) \triangleq \mathsf{E}\{(x(t_f) - x_{t_f})^T P_{t_f}(x(t_f) - x_{t_f}) + \int_t^{t_f}[x^T(s)P(s)x(s)ds + u^T(s,x(s))Q(s)u(s,x(s))]ds\}, \tag{2}$$

where superscript $T$ denotes the transpose, the cost function $J(t,x)$ is the total cost expected to accumulate if the system is initialized in state $x$ at time $t$ and controlled until the final time $t_f$ according to the control law $u = u(t,x)$, $\mathsf{E}\{\cdot\}$ is the expectation with respect to $\sigma$-algebra $\mathfrak{F}$, $x_{t_f}$ is the target position, $P_{t_f}$ and $P(s)$ are symmetric non-negative matrices, $Q(s)$ is a symmetric positive matrix. The quadratic criterion (2) is discussed in our paper to reduce complexity in the optimal control regulator algorithm and reflect physiological characteristics in the biomechanical model.

The optimal control problem is to find the optimal feedback control $u^*(t)$, $0 \le t \le t_f$, that minimizes criterion $J(0, x_0)$, along with the trajectory $x^*(t), 0 \le t \le t_f$, Note that it is difficult to find a globally-optimal control regulator in such the nonlinear stochastic dynamic system. Instead, we try to seek locally-optimal control laws: we will build a Linear-Quadratic-Gaussian (LQG) approximation to our original non-linear time-delay system. Then we will design an quasi-optimal control regulator in section 3.

## 2.2. Local LQG approximation

Throughout the paper, there is an equidistant grid on the time interval $[0, t_f]$ with mesh (sample period) $\Delta t = \frac{t_f}{K} = \frac{\tau}{l} > 0$, for two positive integers $K, l$. The interval $[t_k, t_{k+1}) = [k\Delta t, (k+1)\Delta t)$ is the sampling interval. Assumed that the input is a piecewise constant over the sampling interval, i.e. the zero-order holds assumption holds true: $u(t) = u(k\Delta t) \equiv u_k = constant$, for $k\Delta t \leq t < (k+1)\Delta t$, and similarly for all other time-vary quantities.

Method starts with an open-loop control initial guess $\bar{u}(t)$, and the corresponding "zero-noise" trajectory $\bar{x}(t)$, obtained by applying $\bar{u}(t)$ to the deterministic system $\dot{x}(t) = f(x(t), u(t), u(t-\tau))$ with the initial state $x(0) = x_0$ and initial control input $u(t)=0$ for $t \in [-\tau, 0)$. It can also be done by Euler integration $\bar{x}_{k+1} = \bar{x}_k + f(\bar{x}_k, \bar{u}_k, \bar{u}_{k-l})\Delta t$.

By linearizing the state equation (1) around $\bar{x}, \bar{u}$, we have the discrete-time linear equations, which describe by the state and control deviations $\delta x_k = x_k - \bar{x}_k$, $\delta u_k = u_k - \bar{u}_k$. Written in terms of these deviations, the discrete-time linear state equations approximation to our original non-linear state equation become

$$\delta x_{k+1} = A_k \delta x_k + B_k^0 \delta u_k + B_k^1 \delta u_{k-l} + C_k(\delta u_k, \delta u_{k-l})\xi_k$$

$$k=0, 1, 2, \ldots, K-1. \qquad (3)$$

with $\delta x_0 = 0$, $\delta u_{-1} = \delta u_{-2} = \cdots = \delta u_{-l} = 0$, and the coefficients are

$$A_k = I_n + \Delta t \cdot \partial f / \partial x, \qquad B_k^i = \Delta t \cdot \partial f / \partial u_i, \qquad i = 0, ,$$

$$C_k(\delta u_k, \delta u_{k-l}) = [c_{1k} + C_{1k}^0 \delta u_k + C_{1k}^1 \delta u_{k-l}, \cdots, c_{pk} + C_{pk}^0 \delta u_k + C_{pk}^1 \delta u_{k-l}],$$

$$c_{j,k} = \sqrt{\Delta t} \cdot F^{[j]}, \qquad C_{j,k}^i = \Delta t \cdot \partial F^{[j]} / \partial u_i, \qquad i = 0, , \quad j = 1, \cdots, p$$

$$f = f(\bar{x}_k, \bar{u}_k, \bar{u}_{k-l}), \quad \partial f/\partial x = \frac{\partial f(\bar{x}_k, \bar{u}_k, \bar{u}_{k-l})}{\partial x}, \quad F = F(\bar{u}_k, \bar{u}_{k-l}),$$

$$\partial f/\partial u_0 = \frac{\partial f(\bar{x}_k, \bar{u}_k, \bar{u}_{k-l})}{\partial u_k}, \quad \partial f/\partial u_1 = \frac{\partial f(\bar{x}_k, \bar{u}_k, \bar{u}_{k-l})}{\partial u_{k-l}}$$

$$\partial F^{[i]}/\partial u_0 = \frac{\partial F^{[i]}(\bar{u}_k, \bar{u}_{k-l})}{\partial u_k}, \quad \partial F^{[i]}/\partial u_1 = \frac{\partial F^{[i]}(\bar{u}_k, \bar{u}_{k-l})}{\partial u_{k-l}},$$

$F^{[i]}$ denotes the $i^{th}$ column of $F$, $i=1, \ldots, p$. The noise $\xi_k \sim N(0, I_p)$, $k=0,1,\cdots,K-1$ are independent each other. The $\sqrt{\Delta t}$ term appears because the covariance of Brownian motion grows linearly with time. Since we are interested in multiplicative noises in the control signal, $F$ is only linearized with respect to $\delta u_k$ and $\delta u_{k-l}$. The $i^{th}$ column of the matrix $C_k(\delta u_k, \delta u_{k-l})$ is $c_{ik} + C_{ik}^0 \delta u_k + C_{ik}^1 \delta u_{k-l}$. Thus the noise conditional covariance with respect to $\xi_0, \xi_1, \cdots, \xi_{k-1}$ is

$$Cov[C_k(\delta u_k, \delta u_{k-l})\xi_k] = \sum_{i=1}^p (c_{ik} + C_{ik}^0 \delta u_k + C_{ik}^1 \delta u_{k-l})(c_{ik} + C_{ik}^0 \delta u_k + C_{ik}^1 \delta u_{k-l})^T$$

Similarly, by quadratizing the cost function (2) around $\bar{x}, \bar{u}$, we have a cost-to-go function in the discrete-time quadratic form, for $k= K, K-1,\ldots, 2, 1, 0$,

$$\cos t_k = \tilde{d}_k + 2\delta x_k^T d_k + \delta x_k^T D_k \delta x_k + 2\delta u_k^T e_k + \delta u_k^T E_k \delta u_k, \tag{4}$$

where $\tilde{d}_k = \Delta t(\bar{x}_k^T P(k\Delta t)\bar{x}_k + \bar{u}_k^T Q(k\Delta t)\bar{u}_k)$, $d_k = \Delta t \cdot P(k\Delta t)\bar{x}_k$,

$D_k = \Delta t \cdot P(k\Delta t)$, $e_k = \Delta t \cdot Q(k\Delta t)\bar{u}_k$, $E_k = \Delta t \cdot Q(k\Delta t)$, $k= K-1,\ldots, 1, 0$.

At the final time step $k = K$, it is

$\tilde{d}_k = (\bar{x}(t_f) - x_{t_f})^T P_{t_f}(\bar{x}(t_f) - x_{t_f})$, $d_k = P_{t_f}(\bar{x}(t_f) - x_{t_f})$, $D_k = P_{t_f}$, $e_k = 0$ and $E_k = 0$. So $E_k$ is a symmetric positive matrix for $k \geq 0$. The discrete-time quadratic criterion approximation to our original one becomes

$$J(k) = \mathsf{E}_k \{\sum_{i=k}^K \cos t_i\}, \tag{5}$$

where $\mathsf{E}_k\{\cdot\}$ is the conditional expectation with respect to $\xi_0, \xi_1, \cdots, \xi_{k-1}$. So the LQG

approximation problem with (3) and (5) is to find the optimal control $\delta u_0, \delta u_1, \cdots, \delta u_{K-1}$, that minimizes criterion $J(0)$ (start from beginning 0, answer: initial time is 0), along with the trajectory $\delta x_1, \delta x_2, \cdots, \delta x_K$, generated by the state equation (3).

3. Designing iLQR with input delay

In this section we focus on the above LQG approximation system. For developing an iterative linear Quadratic Regulator (iLQR) algorithm for time-delay system, we first try to find an optimal control $\delta u_{K-1}^*$ to minimize $J(0)$, when the optimal control laws $\delta u_0^*, \delta u_1^*, \cdots, \delta u_{K-2}^*$ are given. Bellman's dynamic programming tells us that it is of equivalence to find an optimal control $\delta u_{K-1}^*$ to minimize $J(K-1)$ in the LQG system. By substituting the state equation (3) into the criterion (5), it gives a new quadratic form $J(K-1)$ without $\delta u_K$ and $\delta x_K$. Notice that the new $J(K-1)$ also has the quadratic form with the controls $\delta u_0^*, \delta u_1^*, \cdots, \delta u_{K-2}^*, \delta u_{K-1}$ and the states $\delta x_{K-1}$. Based upon Pontryagin's maximum principle, the optimal control law $\delta u_{K-1}^*$ can be written as an explicit linear expression with respect to $\delta u_0^*, \delta u_1^*, \cdots, \delta u_{K-2}^*$ and $\delta x_{K-1}$. For simplicity, the optimal control $\delta u_k^*$ is denoted by $\delta u_k$ later. Then repeating the above procedure backwardly, the optimal feedback control laws of the LQG system is derived in the following theorem.

Theorem 1. Consider the LQG system with the state equations (3) and the quadratic criterion (5), the optimal control laws $\delta u_0, \delta u_1, \cdots, \delta u_{K-1}$ satisfy

$$\delta u_k = I_k + L_k \delta x_k + \sum_{i=0}^{\min\{(K-k),l\}-1} M_k^i \delta u_{k+i-l} \qquad (6)$$

And the criterion (5) to be minimum is equivalent to

$$J(k) = \tilde{s}_k + 2\delta x_k^T s_k + \delta x_k^T S_k \delta x_k + 2 \sum_{i=0}^{\min\{(K-k),l\}-1} \delta u_{k+i-l}^T r_k^i$$
$$+ 2\delta x_k^T \left( \sum_{i=0}^{\min\{(K-k),l\}-1} \tilde{R}_k^i \delta u_{k+i-l} \right) + \sum_{i,j=0}^{\min\{(K-k),l\}-1} (\delta u_{k+i-l})^T R_k^{ij} \delta u_{k+j-l} \qquad (7)$$

for $k = 0, 1, \cdots, K-1$. Their coefficients are

$$H_{k-1} = \begin{cases} E_{k-1} + (B_{k-1}^0)^T S_k B_{k-1}^0 + \sum_{i=1}^{p} (C_{i(k-1)}^0)^T S_k C_{i(k-1)}^0, & K \leq k \leq K-l \\ E_{k-1} + (B_{k-1}^0)^T S_k B_{k-1}^0 + \sum_{i=1}^{p} (C_{i(k-1)}^0)^T S_k C_{i(k-1)}^0 + R_k^{(l-1)(l-1)} + (B_{k-1}^0)^T \tilde{R}_k^{l-1} + (\tilde{R}_k^{l-1})^T B_{k-1}^0, & 0 < k \leq K-l-1 \end{cases}$$

$$I_{k-1} = \begin{cases} -(H_{k-1})^{-1} [e_{k-1} + (B_{k-1}^0)^T s_k + \sum_{i=1}^{p} (C_{i(k-1)}^0)^T S_k c_{i(k-1)}], & K \leq k \leq K-l \\ -(H_{k-1})^{-1} [e_{k-1} + (B_{k-1}^0)^T s_k + \sum_{i=1}^{p} (C_{i(k-1)}^0)^T S_k c_{i(k-1)} + r_k^{l-1}], & 0 < k \leq K-l-1 \end{cases}$$

$$L_{k-1} = \begin{cases} -(H_{k-1})^{-1} (B_{k-1}^0)^T S_k A_{k-1}, & K \leq k \leq K-l \\ -(H_{k-1})^{-1} (S_k B_{k-1}^0 + \tilde{R}_k^{l-1})^T A_{k-1}, & 0 < k \leq K-l-1 \end{cases}$$

$$M_{k-1}^0 = \begin{cases} -(H_{k-1})^{-1} [(B_{k-1}^0)^T S_k B_{k-1}^1 + \sum_{i=1}^{p} (C_{i(k-1)}^0)^T S_k C_{i(k-1)}^1], & K \leq k \leq K-l \\ -(H_{k-1})^{-1} [(B_{k-1}^0)^T S_k B_{k-1}^1 + \sum_{i=1}^{p} (C_{i(k-1)}^0)^T S_k C_{i(k-1)}^1 + (\tilde{R}_k^{l-1})^T B_{k-1}^1], & 0 < k \leq K-l-1 \end{cases}$$

$$M_{k-1}^i = \begin{cases} -(H_{k-1})^{-1} (B_{k-1}^0)^T \tilde{R}_k^{i-1}, & K \leq k \leq K-l, \; i=1,2,\cdots,K-k \\ -(H_{k-1})^{-1} [(B_{k-1}^0)^T \tilde{R}_k^{i-1} + R_k^{(l-1)(i-1)}], & 0 < k \leq K-l-1, \; i=1,2,\cdots,l-1 \end{cases}$$

$$\tilde{s}_{k-1} = \tilde{d}_{k-1} + \tilde{s}_k + \sum_{i=1}^{p} (c_{i(k-1)})^T S_k c_{i(k-1)} - I_{k-1}^T H_{k-1} I_{k-1},$$

$$s_{k-1} = d_{k-1} + (s_k)^T A_{k-1} - I_{k-1}^T H_{k-1} L_{k-1},$$

$$S_{k-1} = D_{k-1} + (A_{k-1})^T S_k A_{k-1} - (L_{k-1})^T H_{k-1} L_{k-1},$$

$$r_{k-1}^0 = (B_{k-1}^1)^T s_k + \sum_{i=1}^{p} (c_{i(k-1)})^T S_k C_{i(k-1)}^1 - (M_{k-1}^0)^T H_{k-1} I_{k-1},$$

$$r_{k-1}^i = r_k^{i-1} - (M_{k-1}^i)^T H_{k-1} I_{k-1}, \quad i=1,2,\cdots,\min\{K-k,l-1\},$$

$$\tilde{R}_{k-1}^0 = (B_{k-1}^1)^T S_k (A_{k-1}) - (M_{k-1}^0)^T H_{k-1} L_{k-1},$$

$$\tilde{R}_{k-1}^i = (A_{k-1})^T \tilde{R}_k^{i-1} - (L_{k-1})^T H_{k-1} M_{k-1}^i, \quad i=1,2,\cdots,\min\{K-k,l-1\},$$

$$R_{k-1}^{00} = (B_{k-1}^1)^T S_k B_{k-1}^1 + \sum_{i=1}^{p}(C_{i(k-1)}^1)^T S_k C_{i(k-1)}^1 - (M_{k-1}^0)^T H_{k-1} M_{k-1}^0,$$

$$(R_{k-1}^{0i})^T = R_{k-1}^{i0} = (\tilde{R}_k^{i-1})^T B_{k-1}^1 - (M_{k-1}^0)^T H_{k-1} M_{k-1}^{i-1}, \quad i=1,2,\cdots,\min\{K-k, l-1\},$$

$$R_{k-1}^{ij} = R_k^{(i-1)(j-1)} - (M_{k-1}^{i-1})^T H_{k-1} M_{k-1}^{j-1}, \quad i,j=1,2,\cdots,\min\{K-k, l-1\}. \tag{8}$$

Proof. We want to show that (6) and (7) hold by induction. First of all, for $k = K-1$, from (5) and (4), we know that

$$J(K-1) = \mathsf{E}_{K-1}\{\sum_{i=K-1}^{K} \cos t_i\} = \cos t_{K-1} + \mathsf{E}_{K-1}\{J(K)\}$$

and

$$J(K) = \tilde{s}_K + 2\delta x_K^T s_K + \delta x_K^T S_K \delta x_K$$

with the coefficients

$$\tilde{s}_K = \tilde{d}_K = (\overline{x}(t_f) - x_{t_f})^T P_{t_f}(\overline{x}(t_f) - x_{t_f}),$$

$$s_K = d_K = P_{t_f}(\overline{x}(t_f) - x_{t_f}), \quad S_K = D_K = P_{t_f}. \tag{9}$$

Substituting the state equation (3) into the criterion $J(K-1)$, it gives

$$J(K-1) = \tilde{d}_{K-1} + 2\delta x_{K-1}^T d_{K-1} + \delta x_{K-1}^T D_{K-1} \delta x_{K-1} + 2\delta u_{K-1}^T e_{K-1} + \delta u_{K-1}^T E_{K-1} \delta u_{K-1}$$

$$+ \mathsf{E}_{K-1}\{\tilde{s}_K + 2(A_{K-1}\delta x_{K-1} + B_{K-1}^0 \delta u_{K-1} + B_{K-1}^1 \delta u_{K-1-l} + C_{K-1}(\delta u_{K-1}, \delta u_{K-1-l})\xi_{K-1})^T s_K$$

$$+ (A_{K-1}\delta x_{K-1} + B_{K-1}^0 \delta u_{K-1} + B_{K-1}^1 \delta u_{K-1-l})^T S_K (A_{K-1}\delta x_{K-1} + B_{K-1}^0 \delta u_{K-1} + B_{K-1}^1 \delta u_{K-1-l})$$

$$+ 2(A_{K-1}\delta x_{K-1} + B_{K-1}^0 \delta u_{K-1} + B_{K-1}^1 \delta u_{K-1-l})^T S_K C_{K-1}(\delta u_{K-1}, \delta u_{K-1-l})\xi_{K-1}$$

$$+ (C_{K-1}(\delta u_{K-1}, \delta u_{K-1-l})\xi_{K-1})^T S_K (C_{K-1}(\delta u_{K-1}, \delta u_{K-1-l})\xi_{K-1})\}$$

Using the fact that $trace(UV) = trace(VU)$ and $\xi_k \sim N(0, I_p)$, we have

$$\mathsf{E}_{K-1}\{(C_{K-1}(\delta u_{K-1}, \delta u_{K-1-l})\xi_{K-1})^T S_K (C_{K-1}(\delta u_{K-1}, \delta u_{K-1-l})\xi_{K-1})\}$$

$$= trace\{\sum_{i=1}^{p}(c_{i(K-1)} + C_{i(K-1)}^0 \delta u_{K-1} + C_{i(K-1)}^1 \delta u_{K-1-l})(c_{i(K-1)} + C_{i(K-1)}^0 \delta u_{K-1} + C_{i(K-1)}^1 \delta u_{K-1-l})^T S_K\}$$

$$= \sum_{i=1}^{p}(c_{i(K-1)} + C_{i(K-1)}^0 \delta u_{K-1} + C_{i(K-1)}^1 \delta u_{K-1-l})^T S_K (c_{i(K-1)} + C_{i(K-1)}^0 \delta u_{K-1} + C_{i(K-1)}^1 \delta u_{K-1-l})$$

Then

$$J(K-1) = \tilde{d}_{K-1} + 2\delta x_{K-1}^T d_{K-1} + \delta x_{K-1}^T D_{K-1} \delta x_{K-1} + 2\delta u_{K-1}^T e_{K-1} + \delta u_{K-1}^T E_{K-1} \delta u_{K-1}$$
$$+ \tilde{s}_K + 2(A_{K-1}\delta x_{K-1} + B_{K-1}^0 \delta u_{K-1} + B_{K-1}^1 \delta u_{K-1-l})^T s_K$$
$$+ (A_{K-1}\delta x_{K-1} + B_{K-1}^0 \delta u_{K-1} + B_{K-1}^1 \delta u_{K-1-l})^T S_K (A_{K-1}\delta x_{K-1} + B_{K-1}^0 \delta u_{K-1} + B_{K-1}^1 \delta u_{K-1-l})$$
$$\sum_{i=1}^p (c_{i(K-1)} + C_{i(K-1)}^0 \delta u_{K-1} + C_{i(K-1)}^1 \delta u_{K-1-l})^T S_K (c_{i(K-1)} + C_{i(K-1)}^0 \delta u_{K-1} + C_{i(K-1)}^1 \delta u_{K-1-l})$$

(10)

$$= (\delta u_{K-1})^T H_{K-1} \delta u_{K-1} + 2(\delta u_{K-1})^T G_{K-1} + g_{K-1}$$

where  $H_{K-1} = E_{K-1} + (B_{K-1}^0)^T S_K B_{K-1}^0 + \sum_{i=1}^p (C_{i(K-1)}^0)^T S_K C_{i(K-1)}^0$,

$$G_{K-1} = -H_{K-1}\{I_{K-1} + L_{K-1}\delta x_{K-1} + M_{K-1}^0 \delta u_{k-1-l}\}$$

$$g_{K-1} = \tilde{d}_{K-1} + 2\delta x_{K-1}^T d_{K-1} + \delta x_{K-1}^T D_{K-1} \delta x_{K-1} + \tilde{s}_K + 2(A_{K-1}\delta x_{K-1} + B_{K-1}^1 \delta u_{K-1-l})^T s_K$$
$$+ (A_{K-1}\delta x_{K-1} + B_{K-1}^1 \delta u_{K-1-l})^T S_K (A_{K-1}\delta x_{K-1} + B_{K-1}^1 \delta u_{K-1-l})$$
$$+ \sum_{i=1}^p (c_{i(K-1)} + C_{i(K-1)}^1 \delta u_{K-1-l})^T S_K (c_{i(K-1)} + C_{i(K-1)}^1 \delta u_{K-1-l})$$

$$I_{K-1} = -(H_{K-1})^{-1}[e_{K-1}^T + (B_{K-1}^0)^T s_K + \sum_{i=1}^p (c_{i(K-1)})^T S_K C_{i(K-1)}^0]$$

$$L_{K-1} = -(H_{K-1})^{-1} \left(B_{K-1}^0\right)^T S_K A_{K-1} \qquad (11)$$

$$M_{K-1}^0 = -(H_{K-1})^{-1}[(B_{K-1}^0)^T S_K B_{K-1}^1 + \sum_{i=1}^p (C_{i(K-1)}^0)^T S_K C_{i(K-1)}^1]$$

and $H_{K-1}$ is a symmetric positive matrix (see Appendix). When the optimal control laws $\delta u_0, \delta u_1, \cdots, \delta u_{K-2}$ are given, the well-known Bellman dynamic programming means that

$$\min_{\delta u_0,\cdots,\delta u_{K-1}} J(0) = \min_{\delta u_0,\cdots,\delta u_{K-2}} \mathsf{E}\{\sum_{i=0}^{K-2} \cos t_i + \min_{\delta u_{K-1}} J(K-1)\}.$$

It implies that we shall find the optimal control $\delta u_{K-1}$ which minimizes $J(K-1)$ in the LQG system. Then based upon the Pontryagin's maximum principle, the optimal control law $\delta u_{K-1}$ can be written to a linear expression as follows

$$\delta u_{K-1} = -(H_{K-1})^{-1} G_{K-1} = I_{K-1} + L_{K-1}\delta x_{K-1} + M_{K-1}^0 \delta u_{K-1-l}, \qquad (12)$$

which implies that (6) holds.

Substituting (12) into (10), it yields

$$J(K-1) = \tilde{s}_{K-1} + 2\delta x_{K-1}^T s_{K-1} + \delta x_{K-1}^T S_{K-1} \delta x_{K-1} + 2\delta u_{K-1-l}^T r_{K-1}^0$$
$$+ 2\delta x_{K-1}^T (\tilde{R}_{K-1}^0 \delta u_{K-1-l}) + (\delta u_{K-1-l})^T R_{K-1}^{00} \delta u_{K-1-l}$$

where $\tilde{s}_{K-1} = \tilde{d}_{K-1} + \tilde{s}_K + \sum_{i=1}^{p}(c_{i(K-1)})^T S_K c_{i(K-1)} - I_{K-1}^T H_{K-1} I_{K-1}$,

$$s_{K-1} = d_{K-1} + (A_{K-1})^T s_K - L_{K-1}^T H_{K-1} I_{K-1},$$

$$S_{K-1} = D_{K-1} + (A_{K-1})^T S_K A_{K-1} - (L_{K-1})^T H_{K-1} L_{K-1},$$

$$r_{K-1}^0 = (B_{K-1}^1)^T s_K + \sum_{i=1}^{p}(C_{i(K-1)}^1)^T S_k c_{i(K-1)} - (M_{K-1}^0)^T H_{K-1} I_{K-1},$$

$$\tilde{R}_{K-1}^0 = (A_{K-1})^T S_K (B_{K-1}^1) - (L_{K-1})^T H_{K-1} M_{K-1}^0,$$

$$R_{K-1}^{00} = (B_{K-1}^1)^T S_K B_{K-1}^1 + \sum_{i=1}^{p}(C_{i(K-1)}^1)^T S_K C_{i(K-1)}^1 - (M_{K-1}^0)^T H_{K-1} M_{K-1}^0. \quad (13)$$

So we prove that (6) and (7) hold.

Now assume that (6) and (7) hold for the time step $k$. We want to prove that those hold for the time step $k-1$ in two cases.

Case 1: When $k \in (K-l-1, K-1]$, we know that (7) is

$$J(k) = \tilde{s}_k + 2\delta x_k^T s_k + \delta x_k^T S_k \delta x_k + 2\sum_{i=0}^{K-1-k} \delta u_{k+i-l}^T r_k^i$$
$$+ 2\delta x_k^T (\sum_{i=0}^{K-1-k} \tilde{R}_k^i \delta u_{k+i-l}) + \sum_{i,j=0}^{K-1-k} (\delta u_{k+i-l})^T R_k^{ij} \delta u_{k+j-l}$$

Substituting (3) into the cost function $J(k-1) = \mathsf{E}_{k-1}\{\text{cost}_{k-1} + J(k)\}$, by a similar way, it gives

$$J(k-1) = (\delta u_{k-1})^T H_{k-1} \delta u_{k-1} + 2(\delta u_{k-1})^T G_{k-1} + g_{k-1} \quad (14)$$

where

$$H_{k-1} = E_{k-1} + (B_{k-1}^0)^T S_k B_{k-1}^0 + \sum_{i=1}^{p}(C_{i(k-1)}^0)^T S_k C_{i(k-1)}^0$$

$$G_{k-1} = -H_{k-1}\{I_{k-1} + L_{k-1}\delta x_{k-1} + M_{k-1}^0 \delta u_{k-1-l} + \sum_{i=1}^{K-k} M_{k-1}^i \delta u_{k-1+i-l}\}$$

$$g_{k-1} = \tilde{d}_{k-1} + 2\delta x_{k-1}^T d_{k-1} + \delta x_{k-1}^T D_{k-1} \delta x_{k-1} + \tilde{s}_k + 2(A_{k-1}\delta x_{k-1} + B_{k-1}^1 \delta u_{k-1-l})^T s_k$$

$$+ (A_{k-1}\delta x_{k-1} + B_{k-1}^1 \delta u_{k-1-l})^T S_k (A_{k-1}\delta x_{k-1} + B_{k-1}^1 \delta u_{k-1-l})$$

$$+ \sum_{i=1}^{p}(c_{i(k-1)} + C_{i(k-1)}^1 \delta u_{k-1-l})^T S_k (c_{i(k-1)} + C_{i(k-1)}^1 \delta u_{k-1-l}) + 2\sum_{i=0}^{K-k-1} \delta u_{k+i-l}^T r_k^i$$

$$+ 2(A_{k-1}\delta x_{k-1} + B_{k-1}^1 \delta u_{k-1-l})^T (\sum_{i=0}^{K-k-1} \tilde{R}_k^i \delta u_{k+i-l}) + \sum_{i,j=0}^{K-k-1} (\delta u_{k+i-l})^T R_k^{ij} \delta u_{k+j-l}$$

$$I_{k-1} = -(H_{k-1})^{-1}[e_{k-1} + (B_{k-1}^0)^T s_k + \sum_{i=1}^{p}(C_{i(k-1)}^0)^T S_k c_{i(k-1)}]$$

$$L_{k-1} = -(H_{k-1})^{-1}(B_{k-1}^0)^T S_k A_{k-1} \quad (15)$$

$$M_{k-1}^0 = -(H_{k-1})^{-1}[(B_{k-1}^0)^T S_k B_{k-1}^1 + \sum_{i=1}^{p}(C_{i(k-1)}^0)^T S_k C_{i(k-1)}^1]$$

$$M_{k-1}^i = -(H_{k-1})^{-1}(B_{k-1}^0)^T \tilde{R}_k^{i-1} \quad, i = 1, 2, \cdots, K-k$$

and $H_{k-1}$ is a positive matrix (see Appendix). Pontryagin's maximum principle tells us that the feedback optimal control law

$$\delta u_{k-1} = -(H_{k-1})^{-1} G_{k-1} = I_{k-1} + L_{k-1}\delta x_{k-1} + \sum_{i=0}^{K-k} M_{k-1}^i \delta u_{k-1+i-l} \quad (16)$$

Substituting into the cost function (14), we have

$$J(k-1) = \tilde{s}_{k-1} + 2\delta x_{k-1}^T s_{k-1} + \delta x_{k-1}^T S_{k-1} \delta x_{k-1} + 2\sum_{i=0}^{K-k} \delta u_{k-1+i-l}^T r_{k-1}^i$$

$$+ 2\delta x_{k-1}^T (\sum_{i=0}^{K-k} \tilde{R}_{k-1}^i \delta u_{k-1+i-l}) + \sum_{i,j=0}^{K-k} (\delta u_{k-1+i-l})^T R_k^{ij} \delta u_{k-1+j-l} \quad (17)$$

where $\tilde{s}_{k-1} = \tilde{d}_{k-1} + \tilde{s}_k + \sum_{i=1}^{p}(c_{i(k-1)})^T S_k c_{i(k-1)} - I_{k-1}^T H_{k-1} I_{k-1}$,

$s_{k-1} = d_{k-1} + (A_{k-1})^T s_k - L_{k-1}^T H_{k-1} I_{k-1}$,

$S_{k-1} = D_{k-1} + (A_{k-1})^T S_k A_{k-1} - (L_{k-1})^T H_{k-1} L_{k-1}$,

$r_{k-1}^0 = (B_{k-1}^1)^T s_k + \sum_{i=1}^{p}(C_{i(k-1)}^1)^T S_k c_{i(k-1)} - (M_{k-1}^0)^T H_{k-1} I_{k-1}$,

$$r_{k-1}^i = r_k^{i-1} - (M_{k-1}^i)^T H_{k-1} I_{k-1}, \quad i=1,2,\cdots,K-k,$$

$$\tilde{R}_{k-1}^0 = (A_{k-1})^T S_k (B_{k-1}^1) - (L_{k-1})^T H_{k-1} M_{k-1}^0,$$

$$\tilde{R}_{k-1}^i = (A_{k-1})^T \tilde{R}_k^{i-1} - (L_{k-1})^T H_{k-1} M_{k-1}^i, \quad i=1,2,\cdots,K-k,$$

$$R_{k-1}^{00} = (B_{k-1}^1)^T S_k B_{k-1}^1 + \sum_{i=1}^{p} (C_{i(k-1)}^1)^T S_k C_{i(k-1)}^1 - (M_{k-1}^0)^T H_{k-1} M_{k-1}^0,$$

$$(R_{k-1}^{0i})^T = R_{k-1}^{i0} = (\tilde{R}_k^{i-1})^T B_{k-1}^1 - (M_{k-1}^0)^T H_{k-1} M_{k-1}^{i-1}, \quad i=1,2,\cdots,K-k,$$

$$R_{k-1}^{ij} = R_k^{(i-1)(j-1)} - (M_{k-1}^{i-1})^T H_{k-1} M_{k-1}^{j-1}, \quad i,j=1,2,\cdots,K-k. \tag{18}$$

Case 2: When $k \in (0, K-l-1]$, Similarity, the cost function (7) can be written as follows

$$J(k) = \tilde{s}_k + 2\delta x_k^T s_k + \delta x_k^T S_k \delta x_k + 2\sum_{i=0}^{l-1} \delta u_{k+i-l}^T r_k^i$$
$$+ 2\delta x_k^T (\sum_{i=0}^{l-1} \tilde{R}_k^i \delta u_{k+i-l}) + \sum_{i,j=0}^{l-1} (\delta u_{k+i-l})^T R_k^{ij} \delta u_{k+j-l}$$

Substituting (3) into the cost function $J(k-1) = \mathsf{E}_{k-1}\{\text{cost}_{k-1} + J(k)\}$, by a similar way, we also have the cost function (14) with the coefficients as follows

$$H_{k-1} = E_{k-1} + (B_{k-1}^0)^T S_k B_{k-1}^0 + \sum_{i=1}^{p} (C_{i(k-1)}^0)^T S_k C_{i(k-1)}^0 + R_k^{(l-1)(l-1)} + (B_{k-1}^0)^T \tilde{R}_k^{l-1} + (\tilde{R}_k^{l-1})^T B_{k-1}^0$$

$$G_{k-1} = -H_{k-1}\{I_{k-1} + L_{k-1}\delta x_{k-1} + M_{k-1}^0 \delta u_{k-1-l} + \sum_{i=1}^{l-1} M_{k-1}^i \delta u_{k-1+i-l}\}$$

$$g_{k-1} = \tilde{d}_{k-1} + 2\delta x_{k-1}^T d_{k-1} + \delta x_{k-1}^T D_{k-1}\delta x_{k-1} + \tilde{s}_k + 2(A_{k-1}\delta x_{k-1} + B_{k-1}^1 \delta u_{k-1-l})^T s_k$$
$$+ (A_{k-1}\delta x_{k-1} + B_{k-1}^1 \delta u_{k-1-l})^T S_k (A_{k-1}\delta x_{k-1} + B_{k-1}^1 \delta u_{k-1-l})$$
$$+ \sum_{i=1}^{p} (c_{i(k-1)} + C_{i(k-1)}^1 \delta u_{k-1-l})^T S_k (c_{i(k-1)} + C_{i(k-1)}^1 \delta u_{k-1-l}) + 2\sum_{i=0}^{l-2} \delta u_{k+i-l}^T r_k^i$$
$$+ 2(A_{k-1}\delta x_{k-1} + B_{k-1}^1 \delta u_{k-1-l})^T (\sum_{i=0}^{l-2} \tilde{R}_k^i \delta u_{k+i-l}) + \sum_{i,j=0}^{l-2} (\delta u_{k+i-l})^T R_k^{ij} \delta u_{k+j-l}$$

$$I_{k-1} = -(H_{k-1})^{-1}[e_{k-1} + (B_{k-1}^0)^T s_k + \sum_{i=1}^{p}(C_{i(k-1)}^0)^T S_k c_{i(k-1)} + r_k^{l-1}]$$

$$L_{k-1} = -(H_{k-1})^{-1}(S_k B_{k-1}^0 + \tilde{R}_k^{l-1})^T A_{k-1}$$

$$M_{k-1}^0 = -(H_{k-1})^{-1}[(B_{k-1}^0)^T S_k B_{k-1}^1 + \sum_{i=1}^{p}(C_{i(k-1)}^0)^T S_k C_{i(k-1)}^1 + (\tilde{R}_k^{l-1})^T B_{k-1}^1]$$

$$M_{k-1}^i = -(H_{k-1})^{-1}[(B_{k-1}^0)^T \tilde{R}_k^{i-1} + R_k^{(l-1)(i-1)}] \quad, i=1,2,\cdots,l-1$$

Based upon Pontryagin's maximum principle, the feedback optimal control law can be written as follows:

$$\delta u_{k-1} = -(H_{k-1})^{-1} G_{k-1} = I_{k-1} + L_{k-1}\delta x_{k-1} + \sum_{i=0}^{l-1} M_{k-1}^i \delta u_{k-1+i-l} \tag{19}$$

Substituting into the cost function (14), we have

$$\begin{aligned}
J(k-1) = \tilde{s}_{k-1} + 2\delta x_{k-1}^T s_{k-1} + \delta x_{k-1}^T S_{k-1} \delta x_{k-1} + 2\sum_{i=0}^{l-1}\delta u_{k-1+i-l}^T r_{k-1}^i \\
+ 2\delta x_{k-1}^T (\sum_{i=0}^{l-1}\tilde{R}_{k-1}^i \delta u_{k-1+i-l}) + \sum_{i,j=0}^{l-1}(\delta u_{k-1+i-l})^T R_k^{ij} \delta u_{k-1+j-l}
\end{aligned} \tag{20}$$

where the coefficients $\tilde{s}_{k-1}, s_{k-1}, S_{k-1}, r_{k-1}^i, R_{k-1}^i, \tilde{R}_{k-1}^i, R_{k-1}^{ij}$, $i,j=0,1,\cdots,l-1$ are the same as those in (8).

So (6) and (7) hold by induction. It is complete of the proof.

Appendix

Now we try to prove the matrix $H_k$ defined in (8) is positive. We should switch to a matrix form of the criterion (7) to better represent our solution. So (7) is rewritten as a uniform matrix form as follows:

$$J(k) = (\delta x_k^T \quad \delta u_{k-1}^T \quad \cdots \quad \delta u_{k-l}^T \quad 1)\Gamma_k \begin{pmatrix} \delta x_k \\ \delta u_{k-1} \\ \vdots \\ \delta u_{k-l} \\ 1 \end{pmatrix}, \tag{21}$$

where the block matrix

$$\Gamma_k = \begin{pmatrix} S_k & \tilde{R}_k^{l-1} & \cdots & \tilde{R}_k^0 & s_k \\ (\tilde{R}_k^{l-1})^T & R_k^{(l-1)(l-1)} & \cdots & R_k^{(l-1)0} & r_k^{l-1} \\ \vdots & \vdots & & \vdots & \vdots \\ (\tilde{R}_k^0)^T & (R_k^{(l-1)0})^T & \cdots & R_k^{00} & r_k^0 \\ s_k^T & (r_k^{l-1})^T & \cdots & (r_k^0)^T & \tilde{s}_k \end{pmatrix},$$

the matrices $\tilde{s}_{k-1}, s_{k-1}, S_{k-1}, r_{k-1}^i, R_{k-1}^i, \tilde{R}_{k-1}^i, R_{k-1}^{ij}$, $i,j = 0,1,\cdots,l-1$ are denoted by (8) and $\tilde{R}_k^i = 0, R_k^{ij} = 0, r_k^i = 0$, $i,j = K-k, K-k+1, \cdots, l-1$ when $k \in (K-l-1, K-1]$.

Then by a similar way like that in the proof of Theorem 1, we will prove that $\Gamma_k$ is a non-negative definite symmetric matrix and $H_k$ is a positive definite symmetric matrix by induction.

According to the proof of Theorem 1, (10) can be rewritten as follows:

$$J(K-1) = (\delta x_{K-1}^T \quad \delta u_{K-1}^T \quad 1) \begin{pmatrix} D_{K-1} & 0 & d_{K-1} \\ 0 & E_{K-1} & e_{K-1} \\ d_{K-1}^T & e_{K-1}^T & \tilde{d}_{K-1} \end{pmatrix} \begin{pmatrix} \delta x_{K-1} \\ \delta u_{K-1} \\ 1 \end{pmatrix}$$

$$+ (\delta x_{K-1}^T \quad \delta u_{K-1}^T \quad \delta u_{K-1-l}^T \quad 1) \begin{pmatrix} A_{K-1}^T & 0 \\ (B_{K-1}^0)^T & 0 \\ (B_{K-1}^1)^T & 0 \\ 0 & 1 \end{pmatrix} \begin{pmatrix} S_K & s_K \\ s_K^T & \tilde{s}_K \end{pmatrix} \begin{pmatrix} A_{K-1} & B_{K-1}^0 & B_{K-1}^1 & 0 \\ 0 & 0 & 0 & 1 \end{pmatrix} \begin{pmatrix} \delta x_{K-1} \\ \delta u_{K-1} \\ \delta u_{K-1-l} \\ 1 \end{pmatrix}$$

$$+ \sum_{i=1}^{p} (\delta u_{K-1}^T \quad \delta u_{K-1-l}^T \quad 1) \begin{pmatrix} (C_{i(K-1)}^0)^T \\ (C_{i(K-1)}^1)^T \\ c_{i(K-1)}^T \end{pmatrix} S_K \begin{pmatrix} C_{i(K-1)}^0 & C_{i(K-1)}^1 & c_{i(K-1)} \end{pmatrix} \begin{pmatrix} \delta u_{K-1} \\ \delta u_{K-1-l} \\ 1 \end{pmatrix}$$

$$= (\delta x_{K-1}^T \quad \delta u_{K-1}^T \quad \delta u_{K-1-l}^T \quad 1) \Phi_{K-1} \begin{pmatrix} \delta x_{K-1} \\ \delta u_{K-1} \\ \delta u_{K-1-l} \\ 1 \end{pmatrix}$$

where

$$\Phi_{K-1} = \begin{pmatrix} I & 0 & 0 \\ 0 & I & 0 \\ 0 & 0 & 0 \\ 0 & 0 & 1 \end{pmatrix} \begin{pmatrix} D_{K-1} & 0 & d_{K-1} \\ 0 & E_{K-1} & e_{K-1} \\ d_{K-1}^T & e_{K-1}^T & \tilde{d}_{K-1} \end{pmatrix} \begin{pmatrix} I & 0 & 0 \\ 0 & I & 0 \\ 0 & 0 & 0 \\ 0 & 0 & 1 \end{pmatrix}$$

$$+ \begin{pmatrix} A_{K-1}^T & 0 \\ (B_{K-1}^0)^T & 0 \\ (B_{K-1}^1)^T & 0 \\ 0 & 1 \end{pmatrix} \begin{pmatrix} S_K & s_K \\ s_K^T & \tilde{s}_K \end{pmatrix} \begin{pmatrix} A_{K-1} & B_{K-1}^0 & B_{K-1}^1 & 0 \\ 0 & 0 & 0 & 1 \end{pmatrix}$$

$$+ \begin{pmatrix} 0 & 0 & 0 \\ I & 0 & 0 \\ 0 & I & 0 \\ 0 & 0 & 1 \end{pmatrix} \sum_{i=1}^{p} \begin{pmatrix} (C_{i(K-1)}^0)^T \\ (C_{i(K-1)}^1)^T \\ c_{i(K-1)}^T \end{pmatrix} S_K \begin{pmatrix} C_{i(K-1)}^0 & C_{i(K-1)}^1 & c_{i(K-1)} \end{pmatrix} \begin{pmatrix} 0 & 0 & 0 \\ I & 0 & 0 \\ 0 & I & 0 \\ 0 & 0 & 1 \end{pmatrix}$$

$$= \begin{pmatrix} D_{K-1} & 0 & 0 & d_{K-1} \\ 0 & E_{K-1} & 0 & e_{K-1} \\ 0 & 0 & 0 & 0 \\ d_{K-1}^T & e_{K-1}^T & 0 & \tilde{d}_{K-1} \end{pmatrix}$$

$$+ \begin{pmatrix} A_{K-1}^T S_K A_{K-1} & A_{K-1}^T S_K B_{K-1}^0 & A_{K-1}^T S_K B_{K-1}^1 & A_{K-1}^T s_K \\ (B_{K-1}^0)^T S_K A_{K-1} & (B_{K-1}^0)^T S_K B_{K-1}^0 & (B_{K-1}^0)^T S_K B_{K-1}^1 & (B_{K-1}^0)^T s_K \\ (B_{K-1}^1)^T S_K A_{K-1} & (B_{K-1}^1)^T S_K B_{K-1}^0 & (B_{K-1}^1)^T S_K B_{K-1}^1 & (B_{K-1}^1)^T s_K \\ s_K^T A_{K-1} & s_K^T B_{K-1}^0 & s_K^T B_{K-1}^1 & \tilde{s}_K \end{pmatrix}$$

$$+ \begin{pmatrix} 0 & 0 & 0 & 0 \\ 0 & \sum_{i=1}^p (C_{i(K-1)}^0)^T S_K C_{i(K-1)}^0 & \sum_{i=1}^p (C_{i(K-1)}^0)^T S_K C_{i(K-1)}^1 & \sum_{i=1}^p (C_{i(K-1)}^0)^T S_K c_{i(K-1)} \\ 0 & \sum_{i=1}^p (C_{i(K-1)}^1)^T S_K C_{i(K-1)}^0 & \sum_{i=1}^p (C_{i(K-1)}^1)^T S_K C_{i(K-1)}^1 & \sum_{i=1}^p (C_{i(K-1)}^1)^T S_K c_{i(K-1)} \\ 0 & \sum_{i=1}^p c_{i(K-1)}^T S_K C_{i(K-1)}^0 & \sum_{i=1}^p c_{i(K-1)}^T S_K C_{i(K-1)}^1 & \sum_{i=1}^p c_{i(K-1)}^T S_K c_{i(K-1)} \end{pmatrix}$$

$$= \begin{pmatrix} D_{K-1} + A_{K-1}^T S_K A_{K-1} & A_{K-1}^T S_K B_{K-1}^0 & A_{K-1}^T S_K B_{K-1}^1 & d_{K-1} + A_{K-1}^T s_K \\ (B_{K-1}^0)^T S_K A_{K-1} & E_{K-1} + (B_{K-1}^0)^T S_K B_{K-1}^0 & (B_{K-1}^0)^T S_K B_{K-1}^1 & e_{K-1} + (B_{K-1}^0)^T s_K \\ (B_{K-1}^1)^T S_K A_{K-1} & (B_{K-1}^1)^T S_K B_{K-1}^0 & (B_{K-1}^1)^T S_K B_{K-1}^1 & (B_{K-1}^1)^T s_K \\ d_{K-1}^T + s_K^T A_{K-1} & e_{K-1}^T + s_K^T B_{K-1}^0 & s_K^T B_{K-1}^1 & \tilde{d}_{K-1} + \tilde{s}_K \end{pmatrix}$$

$$+ \begin{pmatrix} 0 & 0 & 0 & 0 \\ 0 & \sum_{i=1}^p (C_{i(K-1)}^0)^T S_K C_{i(K-1)}^0 & \sum_{i=1}^p (C_{i(K-1)}^0)^T S_K C_{i(K-1)}^1 & \sum_{i=1}^p (C_{i(K-1)}^0)^T S_K c_{i(K-1)} \\ 0 & \sum_{i=1}^p (C_{i(K-1)}^1)^T S_K C_{i(K-1)}^0 & \sum_{i=1}^p (C_{i(K-1)}^1)^T S_K C_{i(K-1)}^1 & \sum_{i=1}^p (C_{i(K-1)}^1)^T S_K c_{i(K-1)} \\ 0 & \sum_{i=1}^p c_{i(K-1)}^T S_K C_{i(K-1)}^0 & \sum_{i=1}^p c_{i(K-1)}^T S_K C_{i(K-1)}^1 & \sum_{i=1}^p c_{i(K-1)}^T S_K c_{i(K-1)} \end{pmatrix}$$

$$= \begin{pmatrix} D_{K-1} + A_{K-1}^T S_K A_{K-1} & -H_{K-1} L_{K-1} & A_{K-1}^T S_K B_{K-1}^1 & d_{K-1} + A_{K-1}^T s_K \\ -(L_{K-1})^T H_{K-1} & H_{K-1} & -(M_{K-1}^0)^T H_{K-1} & -(I_{K-1})^T H_{K-1} \\ (B_{K-1}^1)^T S_K A_{K-1} & -H_{K-1} M_{K-1}^0 & (B_{K-1}^1)^T S_K B_{K-1}^1 + \sum_{i=1}^{p}(C_{i(K-1)}^1)^T S_K C_{i(K-1)}^1 & (B_{K-1}^1)^T s_K + \sum_{i=1}^{p}(C_{i(K-1)}^1)^T S_K c_{i(K-1)} \\ d_{K-1}^T + s_K^T A_{K-1} & -H_{K-1} I_{K-1} & s_K^T B_{K-1}^1 + \sum_{i=1}^{p} c_{i(K-1)}^T S_K C_{i(K-1)}^1 & \tilde{d}_{K-1} + \tilde{s}_K + \sum_{i=1}^{p} c_{i(K-1)}^T S_K c_{i(K-1)} \end{pmatrix}$$

and $I$ is a unit matrix with respect to their argument, $H_{K-1}, L_{K-1}, I_{K-1}, M_{K-1}^0$ are defined in (11).

We can obtain that $H_{K-1} = E_{K-1} + (B_{K-1}^0)^T S_K B_{K-1}^0 + \sum_{i=1}^{p}(C_{i(K-1)}^0)^T S_K C_{i(K-1)}^0$ is a positive definite symmetric matrix, because $E_{K-1}$ is a positive definite symmetric matrix and $S_K = P(K\Delta t)$ is a non-negative definite symmetric matrix. By the definition of $D_k, E_k, d_k, e_k, \tilde{d}_k$ in (4) and $S_K, s_K, \tilde{s}_K$ in (9), we know that for $k = 1, 2, \cdots, K-1$,

$$\begin{pmatrix} D_k & 0 & d_k \\ 0 & E_k & e_k \\ d_k^T & e_k^T & \tilde{d}_k \end{pmatrix} = \Delta t \begin{pmatrix} I & 0 \\ 0 & I \\ \bar{x}_k^T & \bar{u}_k^T \end{pmatrix} \begin{pmatrix} P(k\Delta t) & 0 \\ 0 & Q(k\Delta t) \end{pmatrix} \begin{pmatrix} I & 0 & \bar{x}_k \\ 0 & I & \bar{u}_k \end{pmatrix}$$

and

$$\begin{pmatrix} S_K & s_K^T \\ s_K & \tilde{s}_K \end{pmatrix} = \begin{pmatrix} I \\ (\bar{x}(t_f) - x_{t_f})^T \end{pmatrix} P_{t_f} \begin{pmatrix} I & \bar{x}(t_f) - x_{t_f} \end{pmatrix}$$

are two non-negative definite symmetric matrices. Notice that the matrix $S_K$ is non-negative. Then the matrix $\Phi_{K-1}$ is non-negative.

Substituting the feedback optimal control law $\delta u_{K-1} = I_{K-1} + L_{K-1} \delta x_{K-1} + M_{K-1}^0 \delta u_{K-1-l}$ into $J(K-1)$, it gives

$$J(K-1) = (\delta x_{K-1}^T \quad \delta u_{K-1-l}^T \quad 1) \Psi_{K-1} \begin{pmatrix} \delta x_{K-1} \\ \delta u_{K-1-l} \\ 1 \end{pmatrix}.$$

where the block matrix are

$$K = \begin{pmatrix} I & 0 & 0 \\ L_{K-1} & M^0_{K-1} & I_{K-1} \\ 0 & I & 0 \\ 0 & 0 & 1 \end{pmatrix}$$

and

$$\Psi_{K-1} = K^T \Phi_{K-1} K = \begin{pmatrix} S_{K-1} & \tilde{R}^0_{K-1} & s_{K-1} \\ (\tilde{R}^0_{K-1})^T & R^{00}_{K-1} & r^0_{K-1} \\ (s_{K-1})^T & (r^0_{K-1})^T & \tilde{s}_{K-1} \end{pmatrix},$$

with the coefficients $\tilde{s}_{K-1}, s_{K-1}, S_{K-1}, r^0_{K-1}, \tilde{R}^0_{K-1}, R^{00}_{K-1}$ defined in (13).

Let $\Gamma_{K-1} = \Theta^T \Psi_{K-1} \Theta = \Theta^T K^T \Phi_{K-1} K \Theta$, where the block matrix $\Theta = (\theta_{ij})$ with the $n \times n$ unit matrix $\theta_{11} = I$, the $d \times d$ unit matrix $\theta_{2(l+1)} = I$, $\theta_{3(l+2)} = 1$ and the otherwise $\theta_{ij} = 0$. So (7) can be written as (21) with $\tilde{R}^i_{K-1} = 0, R^{ij}_{K-1} = 0, r^i_{K-1} = 0$, $i, j = 1, 2, \cdots, l-1$. We can obtain that the matrix $\Gamma_{K-1}$ is non-negative because the matrix $\Phi_{K-1}$ is non-negative.

Next assume that the matrix $\Gamma_k$ defined in (21) is non-negative. We want to prove that the matrix $\Gamma_{k-1}$ defined in (21) is also non-negative and the matrix $H_k$ is positive.

By a similar way like that in the proof of Theorem 1, substituting (3) into the cost function $J(k-1) = \mathsf{E}_{k-1}\{\text{cost}_{k-1} + J(k)\}$, it yields

$$J(k-1) = (\delta x_{k-1}^T \quad \delta u_{k-1}^T \quad 1) \begin{pmatrix} D_{k-1} & 0 & d_{k-1} \\ 0 & E_{k-1} & e_{k-1} \\ d_{k-1}^T & e_{k-1}^T & \tilde{d}_{k-1} \end{pmatrix} \begin{pmatrix} \delta x_{k-1} \\ \delta u_{k-1} \\ 1 \end{pmatrix}$$

$$+ (\delta x_{k-1}^T \quad \delta u_{k-1}^T \quad \cdots \quad \delta u_{k-l-1}^T \quad 1) \Pi^T \Gamma_k \Pi \begin{pmatrix} \delta x_{k-1} \\ \delta u_{k-1} \\ \vdots \\ \delta u_{k-l-1} \\ 1 \end{pmatrix}$$

$$+ (\delta u_{k-1}^T \quad \delta u_{k-l-1}^T \quad 1) [\sum_{i=1}^{p} \begin{pmatrix} (C_{i(k-1)}^0)^T \\ (C_{i(k-1)}^1)^T \\ c_{i(k-1)}^T \end{pmatrix} S_k (C_{i(k-1)}^0 \quad C_{i(k-1)}^1 \quad c_{i(k-1)})] \begin{pmatrix} \delta u_{k-1} \\ \delta u_{k-l-1} \\ 1 \end{pmatrix}$$

$$= (\delta x_{k-1}^T \quad \delta u_{k-1}^T \quad \cdots \quad \delta u_{k-l-1}^T \quad 1) \Phi_{k-1} \begin{pmatrix} \delta x_{k-1} \\ \delta u_{k-1} \\ \vdots \\ \delta u_{k-l-1} \\ 1 \end{pmatrix}$$

(22)

where

$$\Pi = \begin{pmatrix} A_{k-1} & B_{k-1}^0 & 0 & \cdots & 0 & B_{k-1}^1 & 0 \\ 0 & I & 0 & \cdots & 0 & 0 & 0 \\ \vdots & \vdots & \vdots & & \vdots & \vdots & \vdots \\ 0 & 0 & 0 & \cdots & I & 0 & 0 \\ 0 & 0 & 0 & \cdots & 0 & 0 & 1 \end{pmatrix},$$

the block matrix $\Phi_{k-1} = (\varphi_{ij})$ is symmetric with $\varphi_{ij} = (\varphi_{ji})^T$ for $i > j$ and

$\varphi_{11} = D_{k-1} + A_{k-1}^T S_k A_{k-1}$,

$\varphi_{12} = A_{k-1}^T (S_k B_{k-1}^0 + \tilde{R}_k^{l-1}) = -(H_{k-1} L_{k-1})^T$,

$\varphi_{1j} = A_{k-1}^T \tilde{R}_k^{l+1-j}$, $j = 3, \cdots, l+1$,

$\varphi_{1(l+2)} = A_{k-1}^T S_k B_{k-1}^1$,

$\varphi_{1(l+3)} = d_{k-1} + A_{k-1}^T s_k$,

$\varphi_{22} = E_{k-1} + (B_{k-1}^0)^T S_k B_{k-1}^0 + \sum_{i=1}^{p} (C_{i(k-1)}^0)^T S_k C_{i(k-1)}^0 + R_k^{(l-1)(l-1)} + (B_{k-1}^0)^T \tilde{R}_k^{l-1} + \tilde{R}_k^{l-1} B_{k-1}^0 = H_{k-1}$

$\varphi_{2j} = (B_{k-1}^0)^T \tilde{R}_k^{l+1-j} + R_k^{(l-1)(l+1-j)} = -H_{k-1} M_{k-1}^{l+2-j}$, $j = 3, \cdots, l+1$

$$\varphi_{2(l+2)} = (S_k B_{k-1}^0 + \tilde{R}_k^{l-1})^T B_{k-1}^1 + \sum_{i=1}^{p}(C_{i(k-1)}^0)^T S_k C_{i(k-1)}^1 = -H_{k-1}M_{k-1}^0,$$

$$\varphi_{2(l+3)} = e_{k-1} + (B_{k-1}^0)^T s_k + r_k^{l-1} + \sum_{i=1}^{p}(C_{i(k-1)}^0)^T S_k c_{i(k-1)},$$

$$\varphi_{ij} = R_k^{(l+1-i)(l+1-j)}, \quad i \geq j, \quad i,j = 3,\cdots l+1,$$

$$\varphi_{i(l+2)} = (\tilde{R}_k^{(l+1-i)})^T B_{k-1}^1, \quad i = 3,\cdots l+1,$$

$$\varphi_{i(l+3)} = r_k^{l+1-i}, \quad i = 3,\cdots l+1,$$

$$\varphi_{(l+2)(l+2)} = (B_{k-1}^1)^T S_k B_{k-1}^1 + \sum_{i=1}^{p}(C_{i(k-1)}^1)^T S_k C_{i(k-1)}^1,$$

$$\varphi_{(l+2)(l+3)} = (B_{k-1}^1)^T s_k + \sum_{i=1}^{p}(C_{i(k-1)})^T S_k c_{i(k-1)},$$

$$\varphi_{(l+3)(l+3)} = \tilde{s}_k + \tilde{d}_k + \sum_{i=1}^{p}(c_{i(k-1)})^T S_k c_{i(k-1)}.$$

The block symmetric matrix $\Phi_{k-1} = (\varphi_{ij})$ is non-negative because $\begin{pmatrix} D_{k-1} & 0 & d_{k-1} \\ 0 & E_{k-1} & e_{k-1} \\ d_{k-1}^T & e_{k-1}^T & \tilde{d}_{k-1} \end{pmatrix}$, $\Pi^T \Gamma_k \Pi$ and $S_{k-1}$ are non-negative.

Since $\Pi^T \Gamma_k \Pi$ and $S_{k-1}$ are two non-negative definite symmetric matrices, $(B_{k-1}^0)^T S_k B_{k-1}^0 + \sum_{i=1}^{p}(C_{i(k-1)}^0)^T S_k C_{i(k-1)}^0 + R_k^{(l-1)(l-1)} + (B_{k-1}^0)^T \tilde{R}_k^{l-1} + \tilde{R}_k^{l-1} B_{k-1}^0$ is non-negative. So we also have that $H_{k-1}$ is a positive definite symmetric matrix, because $E_{k-1}$ is a positive definite symmetric matrix.

Notice the relation $\delta u_{k-1} = I_{k-1} + L_{k-1} \delta x_{k-1} + \sum_{i=0}^{l-1} M_{k-1}^i \delta u_{k-1+i-l}$. Let

$$\begin{pmatrix} \delta x_{k-1} \\ \delta u_{k-1} \\ \delta u_{k-1} \\ \vdots \\ \delta u_{k-l-1} \\ 1 \end{pmatrix} = \begin{pmatrix} I & 0 & \cdots & 0 & 0 \\ L_{k-1} & M_{k-1}^{l-1} & \cdots & M_{k-1}^{0} & 0 \\ 0 & I & \cdots & 0 & 0 \\ \vdots & \vdots & & \vdots & \vdots \\ 0 & 0 & 0 & I & 0 \\ 0 & 0 & 0 & 0 & 1 \end{pmatrix} \begin{pmatrix} \delta x_{k-1} \\ \delta u_{k-2} \\ \vdots \\ \delta u_{k-l-1} \\ 1 \end{pmatrix} = K \begin{pmatrix} \delta x_{k-1} \\ \delta u_{k-2} \\ \vdots \\ \delta u_{k-l-1} \\ 1 \end{pmatrix}.$$

Substituting the above equality into (22), we can obtain that $\Gamma_{k-1} = K^T \Phi_{k-1} K$ is non-negative because the matrix $\Phi_{k-1} = (\varphi_{ij})$ is non-negative. The relation $\Gamma_{k-1} = K^T \Phi_{k-1} K$ is equivalence to (8).

In particular, when $k \in (K-l-1, K-1]$, we know that the coefficients $\tilde{R}_k^i = 0, R_k^{ij} = 0, r_k^i = 0$, $i, j = K-k, K-k+1, \cdots, l-1$. Then by carefully calculating from (8), we get that the coefficients $\tilde{R}_{k-1}^i = 0, R_{k-1}^{ij} = 0, r_{k-1}^i = 0$, $i, j = K-k+1, K-k+2, \cdots, l-1$.